\documentclass[11pt,a4paper,reqno,english]{amsart}
\title[]{A dynamical Amrein-Berthier uncertainty principle}
\usepackage{%
    amssymb%
    ,babel%
    ,csquotes%
    ,graphicx%
    ,mathrsfs%
    ,eucal%
    ,xcolor%
    ,enumitem%
    ,geometry%
    ,esint}
\usepackage[%
    pdfauthor={Piero D'Ancona}%
    ,colorlinks=true%
    ,linkcolor=magenta%
    ,citecolor=cyan%
    ,urlcolor=blue%
    ,hyperfootnotes=false%
]{hyperref}

\DeclareMathOperator{\spt}{spt}

\newcommand{\bra}[1]{\langle #1 \rangle}
\newcommand{\one}[1]{\mathbf{1}_{#1}}

\numberwithin{equation}{section}
\newtheorem{theorem}{Theorem}[section]
\newtheorem{corollary}[theorem]{Corollary}
\newtheorem{lemma}[theorem]{Lemma}
\newtheorem{proposition}[theorem]{Proposition}
\theoremstyle{remark}
\newtheorem{remark}[theorem]{Remark}

\theoremstyle{definition}

\newtheorem{example}[theorem]{Example}

\date{\today}

\author[P.~D'Ancona]{Piero D'Ancona}
\address{Piero D'Ancona: 
Dipartimento di Matematica\\
Sapienza Universit\`{a} di Roma\\
Piazzale A.~Moro 2\\
00185 Roma\\
Italy}
\email{dancona@mat.uniroma1.it}
\author[D.~Fiorletta]{Diego Fiorletta}
\address{Diego Fiorletta: 
Dipartimento di Matematica\\
Sapienza Universit\`{a} di Roma\\
Piazzale A.~Moro 2\\
00185 Roma\\
Italy}
\email{diego.fiorletta@uniroma1.it}

\thanks{%
The authors are partially supported by the MIUR PRIN project 2020XB3EFL, ``Hamiltonian and Dispersive PDEs'', 
by the Progetto Ricerca Scientifica 2023
``Long time dynamics of nonlinear systems in non uniform 
environments'' of Sapienza University,and
by the Gruppo Nazionale per l'Analisi Matematica, la Probabilit\`{a} e le loro Applicazioni (GNAMPA), Project CUP E53C23001670001%
}

\begin{document}

\begin{abstract}
  Given a selfadjoint magnetic Schr\"odinger operator 
  \begin{equation*}
    H = ( i \partial + A(x) )^2 + V(x)
  \end{equation*} 
  on $L^{2}(\mathbb{R}^n)$,
  with $V(x)$ strictly subquadratic and $A(x)$ strictly sublinear,
  we prove that the flow $u(t)=e^{-itH}u(0)$ satisfies
  an Amrein--Berthier type inequality
  \begin{equation*}
    \|u(t)\|_{L^{2}}\lesssim_{E,F,T,A,V}
    \|u(0)\|_{L^{2}(E^{c})}
    +
    \|u(T)\|_{L^{2}(F^{c})},
    \qquad
    0\le t\le T
  \end{equation*}
  for all compact sets $E,F \subset \mathbb{R}^{n}$.
  In particular, if both $u(0)$ and $u(T)$
  are compactly supported, then $u$ vanishes identically.
  Under different assumptions on the operator,
  which allow for time--dependent coefficients,
  the result extends to sets $E,F$ of finite measure.
  We also consider a few variants for  Schr\"{o}dinger operators
  with singular coefficients, metaplectic operators,
  and we include applications to control theory.
\end{abstract}

\maketitle



\section{Introduction}\label{sec:intr}

The classical \emph{Heisenberg uncertainty} principle
states that a function $u$ and its Fourier transform
$\widehat{u}(\xi)=\mathcal{F}u=\int e^{-ix \cdot \xi}u(x)dx$
can not be simultaneously localized beyond a certain threshold.
This fact is expressed by the inequality
\begin{equation*}
  \textstyle
  \|xu(x)\|_{L^{2}}\|\xi \widehat{u}(\xi)\|_{L^{2}}
  \ge \frac n2(2\pi)^{\frac n2}\|u\|_{L^{2}}
  \quad\text{for any}\quad 
  u\in H^{1}(\mathbb{R}^{n}),
  \quad
  n\ge1.
\end{equation*}
G.H.~Hardy \cite{Hardy33-a} gave the following
sharp pointwise version:
\begin{equation*}
  \text{if}\ 
  e^{a|x|^{2}}u,\
  e^{\frac{1}{4a}|\xi|^{2}}\widehat{u}\in L^{\infty}
  \ \text{for some}\ 
  a>0
  \quad
   \text{then} 
  \quad
  u=Ce^{-a|x|^{2}},\ 
  \widehat{u}=C'e^{-\frac{1}{4a}|\xi|^{2}}.
\end{equation*}
In particular, if $e^{a|x|^{2}}u$ and 
$e^{b|\xi|^{2}}\widehat{u}$ are bounded and $4ab>1$
then $u$ vanishes identically.
Variants and extensions of Hardy's principle, where
the pointwise conditions on $u,\widehat{u}$ are
replaced by integral conditions, were given in
\cite{Morgan34-a}, 
\cite{Hormander91-a},
\cite{CowlingPrice83-a},
\cite{BonamiDemangeJaming03-a}.
In another direction, Benedicks \cite{Benedicks85-a} proved that
\begin{equation*}
  \text{if $u\in L^{p}$ and}\  
  \spt u,\ \spt \widehat{u}
  \ \text{have finite measure then}
  \quad
  u \equiv0.
\end{equation*}
Note that while Hardy's principle extends to tempered
distributions, the last result requires some restriction;
indeed, it may happen that the supports of
$u,\widehat{u}\in \mathscr{S}'$ have both measure zero,
as shown by the Poisson summation formula
$\mathcal{F}\sum_{j}\delta_{2\pi j}=\sum_{j}\delta_{j}$.
A quantitative version of Benedicks' principle is 
the \emph{Amrein--Berthier inequality}
\cite{AmreinBerthier77-a}:
\begin{equation*}
  \qquad\text{for any $E,F\subset\mathbb{R}^{n}$ of finite measure}
  \quad 
  \|u\|_{L^{2}}\lesssim_{E,F}
  (\|u\|_{L^{2}(E^{c})}+\|\widehat{u}\|_{L^{2}(F^{c})}).
\end{equation*}
It is conjectured that the sharp constant
is of the form $C\exp(C|E|^{1/n}|F|^{1/n})$; this has been
proved in dimension $n=1$ 
\cite{Nazarov93-a}, and only
partial results are available for $n>1$ \cite{Jaming07-a}.

\emph{Dynamical} uncertainty principles stem from the remark
that the Schr\"{o}dinger propagator 
$e^{it \Delta}\phi=
\mathcal{F}^{-1}(e^{-it|\xi|^{2}}\widehat{\phi}(x))$ 
is essentially a Fourier transform: e.g., we can write
\begin{equation*}
  e^{\frac i2\Delta}\phi=
    ce^{-\frac{|x|^{2}}{2i}}
    \mathcal{F}_{y\to x}
      (e^{-\frac{|y|^{2}}{2i}}\phi(y)).
\end{equation*}
Thus all the above principles can be translated immediately 
to this setting, by replacing the couple $u,\widehat{u}$ with
the couple $e^{it \Delta}\phi$, $e^{is \Delta}\phi$ 
at two different times $t\neq s$.
Then a natural question arises, if it is possible to extend
such principles to propagators of
more general evolution equations, like Schr\"{o}dinger equations
with potentials or other dispersive models,
both linear and nonlinear.

The Schr\"{o}dinger equation
\begin{equation}\label{eq:escau}
  iu_{t}+(\Delta+V(t,x))u=0
  \quad\text{on}\quad 
  [0,T]\times \mathbb{R}^{n}
\end{equation}
was investigated from this point of view in a series of papers 
\cite{EscauriazaKenigPonce08-a},
\cite{EscauriazaKenigPonce10-a},
\cite{CowlingEscauriazaKenig10-a}
for bounded potentials $V(t,x)\in L^{\infty}$
satisfying either 
$\int_{0}^{T}\|V(t,x)\|_{L^{\infty}(|x|>R)}dt\to0$ as
$R\to +\infty$, or 
$V=V_{1}(x)+V_{2}(t,x)$ with $V_{1}\in L^{\infty}$ real
and
$|V_{2}(t,x)|\lesssim 
  e^{-(\frac{T}{\alpha t+\beta(T-t)})^{2}|x|^{2}}$
with $\alpha,\beta>0$.
A typical result is that if a solution $u(t,x)$ of 
\eqref{eq:escau} satisfies for some $a,b,T>0$
\begin{equation*}
  e^{a|x|^{2}}u(0,x),\ 
  e^{b|x|^{2}}u(T,x)\in L^{2},
  \qquad
  4abT>1
\end{equation*}
then $u \equiv0$; on the other hand, if
$4abT=1$ one can construct a complex potential
$|V(t,x)|\lesssim \bra{x}^{-2}$ and a nonzero
solution satisfying the above Gaussian bounds.
Corresponding results for the magnetic case
\begin{equation*}
  iu_{t}+(\partial+iA(x))^{2}u+V(t,x)u=0
\end{equation*}
have been proved in 
\cite{BarceloFanelliGutierrez12-a},
\cite{CassanoFanelli15-a}.
In these papers, the assumptions on $V$ are similar to
\cite{EscauriazaKenigPonce08-a}, while the potential 
$A:\mathbb{R}^{n}\to \mathbb{R}^{n}$ is assumed
to be locally integrable along rays emanating from 0,
and such that
$\sum_{j}x_{j}(\partial_{j}A_{k}-\partial_{k}A_{j})\in L^{\infty}$.
See also \cite{Fernandez-BertolinMalinnikova21-a} for an
interesting survey on related problems.

The previous results require $A,V$ to be
bounded, although the Schr\"{o}dinger propagator 
is well defined, and physically relevant,
for more general \emph{subquadratic} potentials $V(t,x)$
and \emph{sublinear} magnetic potentials $A(t,x)$
(see \cite{Yajima91-a}). Thus it is natural to investigate 
the possible extensions of uncertainty principles to
unbounded potentials of this type.
Further evidence is given by the remark that the so--called
\emph{lens transform}
\begin{equation*}
  \textstyle
  u_{le}(t,x)=(\cos t)^{-\frac n2}
  e^{-\frac i4|x|^{2}\tan t}
  u(\tan t,\frac{x}{\cos t}),
  \qquad
  |t|<\frac \pi2
\end{equation*}
takes $u(t,x)=e^{it \Delta}\phi(x)$ 
into a solution of the equation with
a quadratic potential
\begin{equation}\label{eq:harmosc}
  \textstyle
  i \partial_{t}u_{le}+\Delta u_{le}-\frac{|x|^{2}}{4}u_{le}
  =0,
  \qquad
  u_{le}(0,x)=\phi(x).
\end{equation}
Thus it is possible to transfer the previous 
uncertainty principles to the Schr\"{o}dinger
propagator for the harmonic oscillator, and similarly
for propagators with linear magnetic potentials;
see \cite{CassanoFanelli17-a}, 
\cite{Fiorletta23-a}
for results exploring this idea.
Note however that solutions of \eqref{eq:harmosc}
are \emph{periodic} in time, thus a version of the
uncertainty principle is only possible for
\emph{small} times $|T|<\pi$.

In our first result, we consider a selfadjoint Schr\"odinger 
operator on 
$L^{2}(\mathbb{R}^{n})$
\begin{equation} \label{eq:magnHam}
  H = ( i \partial + A(x) )^2 + V(x)
\end{equation}
with \emph{unbounded} potentials $A,V$.
The precise conditions are the following,
where $C_{0}=C_{0}(\mathbb{R}^{n})$ denotes the space of
functions vanishing at infinity. 

\bigskip
\textbf{Assumption (H)}.
$V\in C^{1}(\mathbb{R}^{n};\mathbb{R})$,
$A\in C^{\infty}(\mathbb{R}^{n},\mathbb{R}^{n})$ and
for some $m\ge0$, $\delta\in(0,1)$
\begin{equation}\label{eq:assAV1}
  V\ge-m^{2},
  \qquad
  |A|+\bra{x}^{-\delta}|V|\lesssim \bra{x},
  \qquad
  |\partial A|+\bra{x}^{-\delta}|\partial V|\in C_{0}
\end{equation}
\begin{equation}\label{eq:assA2}
  |\partial^{\alpha}A|\lesssim \bra{x}^{1+\delta(|\alpha|-1)},
  \qquad
  \forall |\alpha|\ge2.
\end{equation}

Thus $V$ and $A$ are allowed to be strictly subquadratic 
and strictly sublinear, respectively. 
Under these conditions, the operator $H$ in
\eqref{eq:magnHam} extends to a selfadjoint operator on 
$L^{2}(\mathbb{R}^{n})$, and standard functional
calculus allows to express the solutions of
$iu_{t}=Hu$ via a continuous group $e^{-itH}$.
Our first result is the following.

\begin{theorem}[Dynamical Amrein--Berthier Inequality]
\label{the:DynAB1}
  Assume $n\ge1$ and the operator $H$ in \eqref{eq:magnHam} 
  satisfies \textbf{(H)}.
  Then for any $T>0$ and any compact sets
  $E$, $F \subset \mathbb{R}^{n}$, 
  there exists a constant $C=C(E,F,T,A,V)$ such that
  any solution of $iu_{t}=Hu$ satisfies
  \begin{equation} \label{eq:AmBerIneq}
    \|u(t)\|_{L^{2}}\le
      C (  \|u(0)\|_{L^{2}(E^{c})}+\|u(T)\|_{L^{2}(F^{c})} )
    \qquad
    \forall t\in \mathbb{R}.
  \end{equation}
\end{theorem}

\begin{remark}[Alternative conditions]\label{rem:altHPpotGlob}
  We list a few variants of Theorem \ref{the:DynAB1}.
  \begin{enumerate}
    \item 
    Theorem \ref{the:DynAB1} holds, more generally, if the potential
    $V$ is of the form $V=V_{1}+V_{2}$, with $V_{1}$
    satisfying \eqref{eq:assAV1} and $V_{2}\in C^{\infty}$ satisfying
    for all $|\alpha|\ge2$
    \begin{equation*}
      V_{2}\ge-m^{2},
      \qquad
      |V_{2}|\lesssim \bra{x}^{2},
      \qquad
      \bra{x}^{-1}|\partial V_{2}|\in C_{0},
      \qquad
      |\partial^{\alpha}V_{2}(x)|\lesssim 
      \bra{x}^{1+\delta(|\alpha|-1)}.
    \end{equation*}

    \item 
    If $A \equiv0$, Theorem \ref{the:DynAB1} holds under
    the assumption $V,\widehat{V}\in L^{1}(\mathbb{R}^{n})$.
    In dimension $n=1$ it is sufficient to assume 
    $V\in L^{1}(\mathbb{R})$.

  \end{enumerate}
\end{remark}

\begin{remark}[]\label{rem:dispEst}
  If the propagator $e^{-itH}$ satisfies a global
  dispersive estimate of the form
  \begin{equation}\label{eq:dispest}
    \|e^{-itH}\phi\|_{L^{p}}\lesssim g(t),
    \qquad
    \lim_{t\to+\infty}g(t)=0
  \end{equation}
  for some $p>2$, then an elementary argument shows that
  \eqref{eq:AmBerIneq} is valid
  for any sets $E,F$ of finite measure provided $T$ is
  large enough with respect to the measures of $E$ and $F$.

  See Section \ref{sub:prfrem} for the proofs of these claims.
\end{remark}

The main drawback of Assumption \textbf{(H)} is the strong
smoothness condition on the coefficients. Using a different
construction of the propagator, we can handle some cases
of singular (and in particular, unbounded) electric potentials:

\begin{theorem}[Dynamical A-B Inequality II]
  \label{the:singPotNonMag}
  Let $H = - \Delta + V(x)$, $n \geq 3$ and let 
  $n^{\star}=\frac{n-1}{n-2}$. 
  Assume that $V \geq - m^2$, and for any $\epsilon > 0$ one has
  $V = V_{1,\epsilon} + V_{2,\epsilon}$, 
  where $\mathcal{F} V_{2,\epsilon}$ is a signed measure 
  of bounded variation, while 
  for some $ \sigma > \frac{1}{n^{\star}}$ and $\gamma > 0$
  \begin{equation} \label{eq:V1Four}
    \lVert \mathcal{F} (\bra{x}^{1 + \sigma} 
      \mathcal{F}^{-1} (\bra{\xi}^{\gamma} 
      \mathcal{F} V_{1,\epsilon})) 
    \rVert_{L^{n^{\star}}(\mathbb{R}^n)}<\epsilon.
  \end{equation}
  Then for any $T>0$ and any compact sets $E$,
  $F \subset \mathbb{R}^n$, solutions of $iu_{t}=Hu$ 
  satisfy \eqref{eq:AmBerIneq}.
\end{theorem}

\begin{remark} \label{rem:Coul}
  To mention an important class of examples, the assumptions of 
  Theorem \ref{the:singPotNonMag} are satisfied by 
  multi--particle Coulomb potentials
  $\sum_{j=1}^{N} \frac{Z_j}{|x - R_j|}$ in dimension $n=3$ 
  and potentials of the form $\sum_{j=1}^{N} \frac{C_{j}}{|x-a_j|^{\frac{n}{n-1}-\delta}}$ in dimension $n \geq 4$. 
\end{remark}

If we assume that the electric potential and the magnetic
field satisfy a suitable symmetry, 
we can greatly enlarge the class of admissible potentials
and allow $V,A$ to depend on time as well.
For instance, we can consider the case of potentials
which are \emph{constant} along a fixed direction $v$
(translational symmetry) or are \emph{rotationally} invariant.
This is the content of the following Theorem \ref{the:magtras}.

\bigskip
\textbf{Assumption (H$_{1}$)}.
$V(t,x)=V_{0}+V_{1}+V_{2}\in 
  C(\mathbb{R}\times \mathbb{R}^{n};\mathbb{R})$
with
$\partial^{\alpha}_{x}V_{0}\in C$ for all 
$\alpha\in \mathbb{N}_{0}$,
$\partial^{\alpha}_{x}V_{0}\in L^{\infty}$ for
all $|\alpha|\ge2$, and for some $p\in(1,\infty]$
\begin{equation*}
  \textstyle
  V_{1}(t,x)\in L^{p}_{t}L^{q}_{x}
  \quad\text{with}\quad 
  \frac 1q=1-\frac{n}{2p}>0,
  \qquad 
  V_{2}\in L^{\infty}_{t}L^{1}_{x},
\end{equation*}
$\partial^{\alpha}_{x}A(t,x)\in 
  C^{1}(\mathbb{R}\times \mathbb{R}^{n};\mathbb{R}^{n})$
for all $|\alpha|\ge0$,
\begin{equation*}
  |\partial^{\alpha}_{x}A(t,x)|+
  |\partial^{\alpha}_{x}\partial_{t}A|+
  \bra{x}^{1+\delta}
  |\partial^{\alpha}_{x}(\partial_{j}A_{k}-\partial_{k}A_{j})|
  \in L^{\infty}
  \quad\text{for all}\quad |\alpha|\ge1.
\end{equation*}

\begin{theorem}[Strong Dynamical A-B Inequality]
  \label{the:magtras}
  Let $H = (i \partial + A(t,x))^2 + V(t,x)$
  be a selfadjoint operator satisfying \textbf{(H$_{1}$)}.
  There exists $\overline{T}\in(0,+\infty]$ such that
  for all $T\in(0,\overline{T}]$ the following holds.
  \begin{enumerate}[label=(\roman*)]
    \item 
    Assume that
    for some direction $v \in S^{n-1}$ and all $j,k=1,\dots,n$
    \begin{equation} \label{eq:compat}
      \partial_{v} V = 
      \partial_{v}\partial_{t}A=
      \partial_{v}( \partial_j A_k - \partial_k A_j ) = 0.
    \end{equation}
    Then, for all $E,F$ of finite measure,
    all solutions of $iu_{t}=Hu$ satisfy \eqref{eq:AmBerIneq}.

    \item 
    Assume $n>1$ and that the operator $H$ commutes with rotations.
    Then, for all $E,F$ of finite measure such that either $E$
    or $F$ does not contain any sphere $\{|x|=R\}$ for $R>0$,
    all solutions of $iu_{t}=Hu$ satisfy \eqref{eq:AmBerIneq}.
  \end{enumerate}
\end{theorem}

\begin{remark} \label{rem:citHP}
  The last equality in condition \eqref{eq:compat} on $A(t,x)$ is implied e.g.~by 
  condition 1.4 in \cite{BarceloFanelliGutierrez12-a} (see also 
  condition 1.7 of Theorem 1.3 in \cite{CassanoFanelli15-a}). 
  Indeed, such condition reads as follows: 
  there exists $v \in S^{n-1}$ such that
  \begin{equation} \label{eq:FanBar}
    \sum_{i=1}^{n} v_i (\partial_i A_l - \partial_l A_i) = 0 ,
    \qquad
    \forall l=1,\ldots,n.
    \end{equation}
    Writing $\partial_v = \sum_{i = 1}^n v_i \partial_i$
    one has:
    \begin{equation*}
    \partial_v (\partial_j A_k - \partial_k A_j) = \sum_{i = 1}^n v_i \partial_i (\partial_j A_k - \partial_k A_j) =  \partial_j  \sum_{i = 1}^n v_i \partial_i A_k - \partial_k \sum_{i=1}^n v_i \partial_i A_j
  \end{equation*}
  which by means of \eqref{eq:FanBar} becomes
  \begin{equation*}
    = \partial_j \sum_{i = 1}^n v_i \partial_k A_i -  \partial_k \sum_{i = 1}^n v_i \partial_j A_i = 0,
    \qquad
    \forall j,k = 1 , \ldots , n.
  \end{equation*}
  On the other hand, to the best of our knowledge, 
  condition (ii) is new.
\end{remark}

\begin{remark}[]\label{rem:heat}
  It is easy to adapt our methods to the heat semigroup
  and prove an estimate of Amrein--Berthier type
  for $u(t)=e^{-tH}u(0)$.
  Indeed, the continuity and uniform boundedness of
  the corresponding kernel is well known from standard results.
  We omit the details.
\end{remark}

Uncertainty principles have many applications to other areas
of PDE theory, notably to control theory.
To show an example of such applications, we deduce
the following observability inequality for 
the propagator $e^{-itH}$.

\begin{theorem} \label{the:Obs1}
  Under the assumptions of Theorem \ref{the:DynAB1} (resp. Theorem \ref{the:singPotNonMag}), for any $T>0$ and for any compact set $E$ we have
  \begin{equation} \label{eq:Obs2}
    \textstyle
    \|u_{0}\|^2_{L^2} \lesssim_{E,T,A,V}
    \int_0^T \| e^{- i t H}u_{0}\|^2_{L^2(E^c)} dt
    \qquad\text{for all}\quad u_{0}\in L^{2}.
  \end{equation}
\end{theorem}

Combining \eqref{eq:Obs2} with the \emph{Hilbert uniqueness method}
(see e.g.~\cite{Zuazua03-a}), one deduces by a standard
procedure the exact controllability property for $e^{-iTH}$.
This means that, given a compact set $E$,
for any $u_0,u_T \in L^2(\mathbb{R}^n)$
we can always find $\nu \in L^2(\mathbb{R}^n)$ such that 
the unique solution of the Cauchy problem
\begin{equation}\label{eq:controleq}
  i \partial_t u(t,x) = H u + \one{E^c} \nu,
  \qquad
  u(0,x)=u_{0}
\end{equation}
reaches the target state at time $T$ i.e. $u(T)=u_T$.

In the final Section \ref{sec:metOp} we explore some applications
of our techniques to \emph{metaplectic operators}, along the lines
of \cite{GrochenigShafkulovska24-a},
\cite{CorderoGiacchiMalinnikova24-a}.

\textbf{Acknowledgments}. The Authors would like to thank Fabio
Nicola and Domenico Monaco for the useful conversations about the topics of this paper.

\textbf{Data Availability Statement}. No new data were created or analysed in this study. Data sharing is not applicable to this article.

\textbf{Conflict Of Interest Statement}. All authors declare that they have no conflicts of interest.

\section{Compactness of the localized propagator}
\label{sec:comp_local_prop}

We denote by $\one{E}$ the characteristic function
of the set $E \subseteq \mathbb{R}^{n}$.
We are interested in the compactness of the localized propagator,
starting from the case of a purely electric potential
$H=-\Delta+V(x)$.

\begin{theorem}[]\label{the:onlyV}
  Let $E,F \subset \mathbb{R}^{n}$ and $H=-\Delta+V(x)$.
  Assume that
  \begin{enumerate}[label=(\roman*)]
    \item either $V(x),\widehat{V}(\xi)\in L^{1}(\mathbb{R}^{n})$
    and $E,F$ have finite measure,
    \item or $n=1$, $V(x)\in L^{1}(\mathbb{R})$ and
    $E,F$ are compact.
  \end{enumerate}
  Then $\one{E}e^{-iTH}\one{F}$ is a compact operator on
  $L^{2}$ for all $T\neq0$.
\end{theorem}

\begin{proof}
  In case (i), the propagator $e^{-iTH}$ is an integral
  operator with bounded kernel $K(T,x,y)$ for every $T\neq0$,
  as proved in Theorem 1.1 from \cite{ZhaoZheng22-a}.
  Hence the kernel $\one{E}(x)K(T,x,y)\one{F}(y)$
  of $\one{E}e^{-iTH}\one{F}$ is in $L^{2}$ and compactness
  follows.

  Consider next case (ii).
  By Theorem B.8.3 from \cite{Simon82-a}, we know that there
  exists a continuous function $K(t,x,y)$ for $t\neq0$ such
  that for all bounded compactly supported $\phi,\psi$ one has
  \begin{equation} \label{eq:weakIntKer}
    \textstyle
    (e^{i T H} \phi , \psi)_{L^2} = 
    \iint_{\mathbb{R}^2} K(T,x,y) \phi(x) \overline{\psi(y)} dx dy,
    \qquad
    T\neq0.
  \end{equation}
  This implies that the operator $M=\one{E} e^{i T H} \one{F}$
  satisfies
  \begin{equation*}
    \textstyle
    (M \phi,\psi)_{L^{2}}=
    \iint K_{E,F}(T,x,y)\phi(x)\overline{\psi(y)}dxdy,
    \qquad
    K_{E,F}=\one{E}(x)K(T,x,y)\one{F}(y)
  \end{equation*}
  where now $K_{E,F}$ is a bounded continuous function
  with compact support, hence in $L^{2}(\mathbb{R}^{2})$.
  Thus the integral operator $\widetilde{M}$ with kernel
  $K_{E,F}$ is a compact operator and
  \begin{equation*}
    ((M-\widetilde{M}) \phi,\psi)=0
  \end{equation*}
  for all bounded compactly supported $\phi,\psi$.
  Since such functions are dense in $L^{2}$ and
  $M-\widetilde{M}$ is a bounded operator, we conclude
  that $M-\widetilde{M}=0$.
\end{proof}

More generally, for magnetic Schr\"{o}dinger operators
\begin{equation*}
  H=(i \partial+A(x))^{2}+V(x)
\end{equation*}
applying the results of \cite{Doi05-a} we obtain the following:

\begin{theorem}\label{the:compact}
  Assume $H$ satisfies \textbf{(H)}, and let 
  $E,F \subset\mathbb{R}^{n}$ be compact sets.
  Then $\one{E}e^{-iTH}\one{F}$ is a compact
  operator on $L^{2}$ for any $T\neq0$.
\end{theorem}

\begin{proof}
  We shall use Theorem 2.5 in \cite{Doi05-a}, concerning
  the smoothness of a Schr\"{o}dinger propagator with
  unbounded potentials.
  In \cite{Doi05-a} a general second order operator
  \begin{equation*}
    L(t)=
    \sum_{j,k=1}^{n}
    (i \partial_{j}+a_{j}(t,x))g^{jk}(x)
    (i \partial_{j}+a_{k}(t,x))
    +V(t,x)
  \end{equation*}
  is considered; our operator $H$ is a special case with
  the choices $a_{j}(t,x)=A_{j}(x)$, $V(t,x)=V(x)$,
  $g^{jk}=\delta_{jk}$.
  We check the assumptions.
  With reference to (H1)--(H5) in \cite{Doi05-a}, we see that
  (H1) and (H2) are trivially satisfied by $g^{jk}=\delta_{jk}$,
  while (H3) follows from \textbf{(H)}.
  Assumption (H5) is satisfied with the choice
  $f_{cv}=1+|x|^{2}$: indeed we have $h_{0}=|\xi|^{2}$
  so that $H_{h_{0}}=2 \xi \cdot \partial_{x}$ and in
  conclusion $H_{h_{0}}^{2}f_{cv}=8|\xi|^{2}$ which
  implies (H5).
  On the other hand, in order to satisfy (H4) it is sufficient
  to assume the following weaker growth condition for $V$,
  \begin{equation}\label{eq:assAV1doi}
    V\ge-m^{2},
    \qquad
    |V|\lesssim \bra{x}^{2},
    \qquad
    \bra{x}^{-1}|\partial V|\in C_{0},
  \end{equation}
  but we must require
  that $V$ is more regular than in Assumption \textbf{(H)},
  namely, we require for the moment that
  \begin{equation}\label{eq:addassV}
    V\in C^{\infty},\qquad
    |\partial^{\alpha}V(x)|\lesssim 
    \bra{x}^{1+\delta(|\alpha|-1)}
    \qquad \forall|\alpha|\ge1.
  \end{equation}
  Thus, for first part of the proof, we assume that
  $V\in C^{\infty}$ satisfies \eqref{eq:assAV1doi} and
  \eqref{eq:addassV}.
  At the end of the proof, by a suitable approximation
  argument, we shall show how to get rid of the
  $C^{\infty}$ assumption and handle $V\in C^{1}$ only.

  In order to apply Theorem 2.5 in \cite{Doi05-a} we choose
  $a(x,\xi)=\chi(x)\in C_{c}^{\infty}(\mathbb{R}^n)$
  an arbitrary test function, so that the corresponding
  Weyl quantization $a^{w}$ is just multiplication by
  $\chi(x)$, and the first coordinate projection
  $\pi(\mathop{\rm supp} a)$ is obviously compact.
  Moreover, the trapped sets $T_{\pm}$ are empty since
  $H_{h_{0}}=2 \xi \cdot \partial_{x}$.
  Finally, the quantities $\mu_{1}(I)=\mu_{2}(I)=0$
  since both $\partial A$ and $\bra{x}^{-1}\partial V$
  vanish at infinity, see \eqref{eq:assAV1}.
  We can thus apply the Theorem with $s=0$, $\rho\ge0$,
  $t_{1}=0$, $t_{2}=T$ and we obtain the estimate
  \begin{equation}\label{eq:doi}
    \|T^{\rho}\bra{x}^{\rho}\chi(x)e^{-iTH}u_{0}\|
      _{L^{\infty}_{T}L^{2}}
    +
    \|T^{\rho}\bra{D}^{\rho}\chi(x)e^{-iTH}u_{0}\|
      _{L^{\infty}_{T}L^{2}}
    \lesssim
    \|\bra{x}^{\rho}u_{0}\|_{L^{2}}.
  \end{equation}
  By compactness of Sobolev embedding,
  from \eqref{eq:doi} we deduce immediately that the
  operator $\chi(x)e^{-iTH}\one{F}$ is compact from
  $L^{2}$ to $L^{2}$ for any test function $\chi$
  and any compact subset $F \subset \mathbb{R}^n$
  and in particular the claim follows.

  It remains to get rid of assumption \eqref{eq:assAV1doi}, \eqref{eq:addassV}.
  If $V$ satisfies only \textbf{(H)}, we define a regularized
  potential $V_{\epsilon}\in C^{\infty}$ as follows.
  Fix a dyadic partition of unity $(\phi_{j})_{j\ge0}$,
  that is to say $\phi_{0}(x)=\chi(x)$,
  $\phi_{j}(x)=\chi(2^{-j}x)-\chi(2^{1-j}x)$
  for a given $\chi\in C_{c}^{\infty}(\mathbb{R}^{n})$
  satisfying $\one{B(0,1)}\le \chi\le \one{B(0,2)}$,
  so that $\sum_{j\ge0}\phi_{j}(x)=1$.
  For $\epsilon>0$ we define then
  \begin{equation*}
    V_{\epsilon}(x)=
    \sum_{j\ge0}\phi_{j}(x)\cdot\rho_{\epsilon_{j}}*V(x),
    \qquad
    \epsilon_{j}=2^{-\delta j}\epsilon
  \end{equation*}
  where $\rho_{\epsilon}(x)=\epsilon^{-n}\rho(x/\epsilon)$
  is the usual Friedrichs mollifier.
  We have
  \begin{equation*}
    \textstyle
    |V(x)-\rho_{\epsilon_{j}}*V(x)|\le
    \int \rho_{\epsilon_{j}}(y)|V(x-y)-V(x)|dy=
    \int \rho_{\epsilon_{j}}(y)|\partial V(\xi_{x,y})||y|dy
  \end{equation*}
  for some point $\xi_{x,y}$ on the segment joining $x,x-y$.
  Since $y\in \spt \rho_{\epsilon_{j}}$ we have
  $|y|\le 2^{-\delta j}\epsilon$, and 
  by assumption \textbf{(H)} we have
  $|\partial V(\xi_{x,y})|\lesssim \langle \xi_{x,y} \rangle^{\delta}
    \lesssim \bra{x}^{\delta}$
  provided $\epsilon\in(0,1)$.
  This implies
  \begin{equation*}
    |V(x)-\rho_{\epsilon_{j}}*V(x)|\lesssim
      2^{-\delta j}\epsilon \bra{x}^{\delta}
  \end{equation*}
  so that, for $x\in\spt \phi_{j}$,
  \begin{equation*}
    \textstyle
    |V(x)-V_{\epsilon}(x)|
    \lesssim
    \epsilon
    2^{-\delta j}\bra{x}^{\delta}
    \lesssim \epsilon.
  \end{equation*}
  Summing up we have
  \begin{equation}\label{eq:VVep}
    \|V-V_{\epsilon}\|_{L^{\infty}}\lesssim \epsilon.
  \end{equation}
  Moreover, we have obviously $V_{\epsilon}\in C^{\infty}$.
  We check that $V_{\epsilon}$ satisfies 
  \eqref{eq:assAV1doi}, \eqref{eq:addassV}.
  First of all we have obviously
  $V_{\epsilon}(x)\ge-m^{2}\sum \phi_{j}=-m^{2}$ and
  \begin{equation}\label{eq:estVep}
    \textstyle
    |\rho_{\epsilon_{j}}*V(x)|\le
    \int \rho_{\epsilon_{j}}(y)|V(x-y)|dy\le
    \sup_{|y|\le 1}|V(x-y)|\lesssim \bra{x}^{2}
  \end{equation}
  so that
  \begin{equation*}
    |V_{\epsilon}(x)|\lesssim \bra{x}^{2}.
  \end{equation*}
  The derivative $\partial V_{\epsilon}$ gives
  two terms: the first one can be estimated like
  \eqref{eq:estVep}
  \begin{equation*}
    \textstyle
    |\sum \phi_{j}\rho_{\epsilon_{j}}*(\partial V)(x)|
    \lesssim
    \bra{x}^{\delta}
  \end{equation*}
  using the assumption $|\partial V(x)|\lesssim \bra{x}^{\delta}$
  from \eqref{eq:assAV1}; the second term gives
  \begin{equation*}
    \textstyle
    |\sum \partial \phi_{j}(x)\rho_{\epsilon_{j}}*V(x)|
    \lesssim
    \sum 2^{-j}\one{C_{j}}(x) \cdot\bra{x}^{1+\delta}
    \lesssim \bra{x}^{\delta},
    \qquad
    C_{j}=B(0,2^{j+1})\setminus B(0,2^{j-1})
  \end{equation*}
  which implies \eqref{eq:assAV1doi} for $V_{\epsilon}$.

  Consider next \eqref{eq:addassV}.
  The derivative $\partial^{\alpha}V_{\epsilon}$
  produces several terms; the worst one is
  \begin{equation*}
    \textstyle
    |\sum \phi_{j}\partial^{\alpha}(\rho_{\epsilon_{j}}*V)|\le
    \sum \phi_{j}|(\partial^{\alpha-e_{k}}\rho_{\epsilon_{j}})
      *(\partial_{k}V)|
  \end{equation*}
  where $e_{k}$ is the $k$--th element of the canonical basis
  of $\mathbb{R}^{n}$, provided $\alpha_{k}\ge1$.
  This gives
  \begin{equation*}
    \textstyle
    \lesssim
    \sum \phi_{j}(x)\epsilon_{j}^{1-|\alpha|}
    \bra{x}^{\delta}=
    \epsilon^{1-|\alpha|}
    \sum \phi_{j}(x)2^{(|\alpha|-1)\delta j}
    \bra{x}^{\delta}\lesssim
    \bra{x}^{\delta|\alpha|}
  \end{equation*}
  which implies \eqref{eq:addassV}.

  If we denote by $H_{\epsilon}$ the
  operator $H$ with $V$ replaced by $V_{\epsilon}$,
  by the first part of the proof we see that
  $\one{E}e^{-iTH_{\epsilon}}\one{F}$ is a compact operator.
  Moreover by Duhamel's principle we can write
  \begin{equation*}
    \textstyle
 \one{E}e^{-iTH}\one{F}=
    \one{E}e^{-iTH_{\epsilon}}\one{F}
    +i\int_{0}^{T}\one{E}e^{-i(T-t)H_{\epsilon}}
      (H-H_{\epsilon})e^{-itH}\one{F}dt
  \end{equation*}
  and the second term at the right is a bounded operator on 
  $L^{2}$, with norm bounded by
  \begin{equation*}
    \lesssim\|H-H_{\epsilon}\|_{L^{2}\to L^{2}}\le
    \|V-V_{\epsilon}\|_{L^{\infty}}\lesssim \epsilon.
  \end{equation*}
  Letting $\epsilon\to0$, we see that 
  $e^{iTH_{\epsilon}}\to e^{iTH}$ in the operator norm,
  hence $\one{E}e^{-iTH}\one{F}$ is the uniform limit
  of the compact operators
  $\one{E}e^{-iTH_{\epsilon}}\one{F}$
  and the proof is concluded.
\end{proof}

\begin{remark}\label{rem:V1V2}
  From the previous proof, we see that conclusion of
  Theorem \ref{the:compact} holds if we assume more
  generally that $V=V_{1}+V_{2}$ with $V_{1}$
  satisfying \eqref{eq:assAV1} and $V_{2}$ satisfying
  \eqref{eq:assAV1doi}, \eqref{eq:addassV}.
\end{remark}

If we restrict to \emph{small times}, we can handle
the more general time dependent operators
\begin{equation*}
  H(t)=(i \partial+A(t,x))^{2}+V(t,x).
\end{equation*}
We define $H(t)$ e.g.~on $C_{c}^{\infty}(\mathbb{R}^{n})$,
and we consider the corresponding Cauchy problem
\begin{equation}\label{eq:CPt}
  iu_{t}+H(t)u=0,
  \qquad
  u(0,x)=\phi(x).
\end{equation}
We denote the solution to \eqref{eq:CPt} by $U(t,0)\phi$.
More generally, we denote by $U(t,s)\phi$ the solution
to the Cauchy problem
\begin{equation}\label{eq:Uts}
  iu_{t}(t,s,x)+H(t)u(t,s,x)=0,
  \qquad
  u(s,s,x)=\phi(x).
\end{equation}

We summarize the relevant results from \cite{Yajima91-a}:

\begin{theorem}\label{the:yaji}
  Assume $A,V$ satisfy conditions \textbf{(H$_{1}$)}
  (see the Introduction).
  Then for all $t\in \mathbb{R}$ the operator $H(t)$ is 
  essentially selfadjoint with domain independent of $t$,
  and for all $\phi\in L^{2}$
  the Cauchy problem \eqref{eq:CPt} has a unique global solution
  $u\in C(\mathbb{R};L^{2})$.

  Moreover, there exists $T>0$ such that the propagator
  $U(t,0)\phi=u(t,x)$ for $t\in[0,T]$
  is an integral operator with bounded kernel
  \begin{equation*}
    \textstyle
    u(t,x)=U(t,0)\phi(x)=t^{-n/2}
    \int K(t,x,y)\phi(y)dy,
    \qquad
    K(t,x,y)\in L^{\infty}([0,T]\times \mathbb{R}^{2n}_{x,y}).
  \end{equation*}
\end{theorem}

\begin{proof}
  This is a consequence of
  Theorems 1--4 and Remark (a) from \cite{Yajima91-a}.
  In particular, it is proved
  that the solution $U(t,s)\phi(x)$ of problem
  \eqref{eq:Uts} is well defined for all $t,s$,
  and the propagator $U(t,s)$
  can be represented for $0<|t-s|\le T$, with $T$ small enough,
  as an oscillatory integral operator
  \begin{equation*}
    \textstyle
    U(t,s)\phi(x)=
    (4\pi i(t-s))^{-n/2}
    \int e^{iS(t,s,x,y)}
    b(t,s,x,y)\phi(y)dy
  \end{equation*}
  (see Theorem 4 in \cite{Yajima91-a}). We have
  $\partial^{\alpha}_{x,y}b(t,s,x,y)\in 
    C^{1}\cap L^{\infty}$
  for all $\alpha\in \mathbb{N}^{2n}_{0}$,
  while the phase $S$ is real valued, $C^{1}$ in all variables,
  $C^{\infty}$ in $(x,y)$, and satisfies
  \begin{equation*}
    \partial_{t}S(t,s,x,y)+(\partial_{x}S(t,s,x,y)-A(t,x))^{2}=0,
  \end{equation*}
  \begin{equation*}
    \textstyle
    \partial^{\alpha}_{x,y}(S(t,s,x,y)-\frac{|x-y|^{2}}{4(t-s)})
    \in L^{\infty}
    \quad\text{for all}\quad |\alpha|\ge2.
  \end{equation*}
  All claims of the above Theorem \ref{the:yaji}
  are easy consequences of these results.
\end{proof}

Consider next the localized operator
\begin{equation*}
  M=\chi(x)U(t,0)\psi(x) 
\end{equation*}
for two bounded functions $\chi,\psi$.
By Theorem \ref{the:yaji}, $M$ is an integral operator with
kernel $\chi(x)K(t,x,y)\psi(y)\in L^{\infty}$.
If we assume that
$E=\spt \chi$ and $F=\spt \psi$ have finite measure,
we have in addition
\begin{equation*}
  \textstyle
  \iint|\chi(x)K(t,x,y)\psi(y)|^{2}dxdy \lesssim
  |E||F|
  \|K\|_{L^{\infty}}^{2}.
\end{equation*}
This implies that $M$ is a Hilbert--Schmidt operator,
and in particular we obtain:

\begin{corollary}\label{cor:compact2}
  Let $\chi,\psi\in L^{\infty}(\mathbb{R}^{n})$
  with supports of finite measure,
  and let $T$ be as in Theorem \ref{the:yaji}.
  Then $\chi(x)U(t,0)\psi(x)$ is a compact operator on $L^{2}$
  for every $t\in(0,T]$.
\end{corollary}

\section{Proof of the main results}\label{sec:proomainres}

\subsection{Proof of Theorem \ref{the:DynAB1}}\label{sec:prfTh1}

Define by spectral calculus
\begin{equation*}
    C(\tau)=\cos(\tau \sqrt{H+m^2}),
    \qquad
    \tau\in \mathbb{R}.
\end{equation*}
We shall need the following elementary property:

\begin{proposition}[Finite Propagation Speed]\label{pro:FPS}
  Assume $\phi\in L^{2}(\mathbb{R}^n)$ is supported in a ball
  $B(x_{0},R)\subset \mathbb{R}^n$. Then
  $C(\tau)\phi$ is supported in
  $B(x_{0},R+|\tau|)$ for all $\tau\in \mathbb{R}$.
\end{proposition}

\begin{proof}
  By a standard approximation argument, we can assume that
  $\phi\in C_{c}^{\infty}$ and 
  $u(t,x)=\cos(t\sqrt{H+m^{2}})\phi\in C([-T,T];H^{1})$.
  Then $u(t,x)$ solves
  \begin{equation*}
    u_{tt}+(H+m^{2})u=0,
    \qquad
    u(0,x)=\phi,
    \quad
    u_{t}(0,x)=0
  \end{equation*}
  and the finite speed of propagation follows by the classical
  local energy estimate on cones of propagation.
\end{proof}

We shall also need a classical result on matrix algebras:

\begin{theorem}[Burnside]\label{the:burn}
  Every family of commuting operators in a finite dimensional
  complex vector space has a common eigenvector.
\end{theorem}

Given the compact sets $E,F$, we know from 
Theorem \ref{the:compact} that the operator
\begin{equation*}
  S=\one{E}e^{iTH}\one{F}e^{-iTH}
\end{equation*}
is a compact operator on $L^{2}$, and we have obviously
\begin{equation}\label{eq:norm1}
  \|S\|_{L^{2}\to L^{2}}\le 1.
\end{equation}
We shall prove that the inequality \eqref{eq:norm1} is strict, 
and Theorem \ref{the:DynAB1} will follow.

Pick a sequence $g_{j}\in L^{2}$
with $\|g_{j}\|_{L^{2}}=1$ such that
$\|Sg_{j}\|_{L^{2}}\to 1$; up to a subsequence, we can assume
that $g_j$ converges weakly to some $ g \in L^{2}$ with $\|g\|_{L^{2}}\le1$.
By compactness of $S$ this implies that $Sg_{j}$ has a
convergent subsequence $Sg_{j_{k}}\to Sg$.
Thus we have $1=\|Sg\|_{L^{2}}\le\|g\|_{L^{2}}\le1$,
so that $\|Sg\|_{L^{2}}=\|g\|_{L^{2}}=1$.
Let $h=e^{iTH}\one{F}e^{-iTH}g$; since
$Sg=\one{E}h$ has norm 1 and $\|h\|_{L^{2}}\le1$,
we see that $h$ must have support contained in $E$,
so that $Sg=h$ has actually norm 1.
It follows that $\one{F}e^{-iTH}g$ has also norm 1
and by the same argument this implies that 
$e^{-iTH}g$ has support in $F$. Summing up, we can write
\begin{equation*}
  Sg=\one{E}e^{iTH}\one{F}e^{-iTH}g=
  e^{iTH}e^{-iTH}g=g
\end{equation*}
and we conclude that $g$ is an eigenfunction of $S$ corresponding
to the eigenvalue 1; moreover, $g$ is supported in $E$ and
$e^{-iTH} g$ is supported in $F$.

Let us consider the family of operators:
\begin{equation}
\label{eqn:opEnl}
S_{t} = \one{\bar{B}_{R+t}(0)} e^{i T H} \one{\bar{B}_{R+t}(0)} e^{-i T H} 
\end{equation}
where $R$ is large enough, i.e. such that $E$,$F \subset B_{R}(0) $, and:
\begin{equation}
\label{eqn:eigSpace}
A^{1}_{t} = \{ f \in L^2 | S_{t} f = f \}
\end{equation}
The $S_t$ are compact, again by means of Theorem \ref{the:compact},
while the set inclusion:
\begin{equation*}
\bar{B}_{R+t_1}(0) \subset \bar{B}_{R+t_2}(0) \text{ if } t_1 < t_2
\end{equation*}
yields the corresponding inclusion of the $A^1_t$:
\begin{equation}
\label{eq:eigInc}
A^1_{t_1} \subset A^1_{t_2} \text{ if } t_1 < t_2
\end{equation}
Thus the map
\begin{equation}\label{eq:dimFunc}
  t \mapsto \operatorname{dim} (A^1_t)
  \quad\text{on}\quad 
  \mathbb{R}_{\geq 0} \rightarrow \mathbb{N} \setminus \{0\}
\end{equation}
is a non-decreasing function. Note that $\operatorname{dim} ( A^1_{t} ) \geq 1$, since we proved the existence of at least an eigenfunction $g$.
Moreover, the compactness of the operators $S_t$ implies that the eigenspaces $A^1_t$ are finite dimensional; this implies that the function \eqref{eq:dimFunc} is piecewise constant. But this and \eqref{eq:eigInc} yield that we can find $\delta$ such that $A^1_{t_1} = A^1_{t_2}$ for each $t_1,t_2 \in (0,\delta)$. Denote by $N$ the dimension of this subspace and consider a basis $f_1 , \ldots , f_N$ of $A^1_{\epsilon}$ for $0 < \epsilon < \delta$, which is also a basis for $A^1_{t}=A^1_{\epsilon}$, $\forall 0 < t < \epsilon$. So $f_1 , \ldots , f_N$ are $1$-eigenfunctions for each operator $S_t$, $0<t<\epsilon$. But this is possible if and only if:
\begin{equation}
\label{suppInc}
\operatorname{supp} ( f_i) \cup \operatorname{supp} (e^{- i T H} f_i) \subset \bigcap_{0 < t < \epsilon} \bar{B}_{R + t} = \bar{B}_{R} \text{ , } \forall i = 1, \ldots , N.
\end{equation}
which in turn means that $f_1 , \ldots , f_N$ is a basis of $A^1_{0}$, and thus $\operatorname{dim}(A^1_0)=N$.
We have thus proved that $\exists \delta > 0$ such that $A^1_{0} = A^1_{t}$, $\forall 0 \leq t < \delta$.

Now, we see that
\begin{equation}
\label{eqn:eigFPS}
\cos(\tau \sqrt{H+m^2}) A^{1}_0 \subset A^{1}_{\tau}=A^1_{0} \text{ , } \forall 0 \leq \tau < \delta
\end{equation}
and thus the cosine operators map $1$-eigenfunctions of $S_0$ into $1$--eigenfunctions of $S_0$. 
Indeed, the inclusion in \eqref{eqn:eigFPS} is easily verified, recalling that $\cos(\tau \sqrt{H+m^2})$ and $e^{- i T H}$ commute and the finite speed of propagation proved in Proposition \ref{pro:FPS}. 

Since the operators $\cos(\tau \sqrt{H+m^2})$ for $0\le \tau<\delta$ form a commuting family of operators acting on the finite dimensional subspace $A^1_{0}$ of $L^2$, by Theorem \ref{the:burn} they share a common eigenfunction $f\in A^1_0$: there exists $\lambda(\tau)\in \mathbb{R}$ such that
\begin{equation} \label{eqn:eigComm}
  \cos(\tau \sqrt{H+m^2}) f = \lambda(\tau) f,
  \qquad
  0\le \tau<\delta.
\end{equation}
In a certain sense, we have built a standing wave solution to the magnetic wave equation, with compact support.

By functional calculus, we see that $\cos(\tau \sqrt{H+m^2}) f = \lambda(\tau)f \in C([0,\delta];L^2(\mathbb{R}^n))$ (this follows from the spectral theorem and dominated convergence e.g.~as in Example 4.12 of \cite{Amrein09-a}), so that $\lambda(\tau) \in C([0,\delta])$. Since $\lambda(0)=1$ we see that $\lambda(\tau)>0$ for $\tau$ close to 0.
We may regard $\lambda(\tau) f(x)$ as a solution to the magnetic wave equation in the sense of distributions:
\begin{equation}
\label{eqn:disWaveEq}
\partial^2_{\tau} \lambda(\tau) f(x) = - \lambda(\tau) (H+m^2) f(x).
\end{equation}
Dividing by $\lambda$ we get
\begin{equation*}
\frac{ \partial^2_{\tau} \lambda(\tau)}{- \lambda(\tau)} f(x) = 
(H+m^2) f(x)
\end{equation*}
which implies that, for some constant $c\in \mathbb{R}$,
\begin{equation*}
\frac{ \partial^2_{\tau} \lambda(\tau)}{\lambda(\tau)} = - c.
\end{equation*}
This proves that $f$ satisfies in the sense of distributions the identity
\begin{equation}
\label{eqn:eigEq}
(H+m^2)f=cf.
\end{equation}
Since $f$ is compactly supported, by unique continuation 
(applying e.g.~Corollary 1.4 from \cite{LuLv17-a} on arbitrarily 
large balls containing the support of $f$) we obtain
that $f$ is identically 0, contradicting Theorem \ref{the:burn}.
We conclude that the norm of the operator $S$ is \emph{strictly}
less than 1.

We can now prove the main inequality \eqref{eq:AmBerIneq}. First of all, the adjoint $S^*$ also satisfies
\begin{equation*}
  \|S^{*}\|_{L^{2}\to L^{2}} = 
  \|e^{i T H} \one{F} e^{-i T H} \one{E} \|_{L^{2}\to L^{2}} <1.
\end{equation*}
Thus we can write
\begin{align*} 
  \|u(0)\|_{L^{2}}  
  &=
  \|e^{-iTH}u(0)\|_{L^{2}}
\\
    &=
    \|e^{-iTH}u(0)\|_{L^{2}(F)}+
    \|e^{-iTH}u(0)\|_{L^{2}(F^c)}
\\ 
    &=
    \|e^{iTH}\one{F}e^{-iTH}u(0)\|_{L^{2}}+
    \|u(T)\|_{L^{2}(F^c)}
\\
  &\le
  \|e^{iTH}\one{F}e^{-iTH}\one{E}u(0)\|_{L^{2}}+
  \|e^{iTH}\one{F}e^{-iTH}\one{E^c}u(0)\|_{L^{2}}+
  \|u(T)\|_{L^{2}(F^c)}
\\
  &\le
  \|S^*u(0)\|_{L^{2}}+
  \|u(0)\|_{L^{2}(E^c)}+  \|u(T)\|_{L^{2}(F^c)}
\\
  &\le\|S^*\|_{L^{2}\to L^{2}}\|u(0)\|_{L^{2}}+
  \|u(0)\|_{L^{2}(E^c)}+  \|u(T)\|_{L^{2}(F^c)}.
\end{align*}
Since $\|S^*\|<1$ and $  \|u(t)\|_{L^{2}}=\|u(0)\|_{L^{2}}$
for all $t$, this implies \eqref{eq:AmBerIneq} with a constant
\begin{equation*}
  C=(1-\|S^*\|_{L^{2}\to L^{2}})^{-1}.
\end{equation*}
This concludes the proof of Theorem \ref{the:DynAB1}.

\subsection{Proof of Remarks \ref{rem:altHPpotGlob},
\ref{rem:dispEst}} \label{sub:prfrem}

\label{rem:variants}
The scheme of Theorem \ref{the:DynAB1} can be easily
adapted to the other cases listed in Remark 
\ref{rem:altHPpotGlob}. It is sufficient to note that:
\begin{itemize}
  \item 
  Compactness of 
  $ \one{E} e^{i T H} \one{F} e^{-i T H}$ is guaranteed 
  by Remark \ref{rem:V1V2} or Theorem \ref{the:onlyV}.

  \item 
  The unique continuation property 
  for solutions of $(H + m^2) \phi = c \phi$ 
  holds in all cases under consideration; in particular,
  in the case $n=1$, $V\in L^{1}$ one can apply
  Theorem 7.1 from \cite{SchechterSimon80-a}.
\end{itemize}

On the other hand, the proof of Remark \ref{rem:dispEst} 
is elementary. Indeed, we can estimate the operator
$S=\one{E}e^{iTH}\one{F}e^{-iTH}$ as follows
\begin{equation*}
  \|\one{E}e^{iTH}\one{F}e^{-iTH}\phi\|_{L^{2}}\le
  |E|^{\frac 12-\frac 1p}\|e^{-iTH}\|_{L^{q}\to L^{p}}
  |F|^{\frac 1q-\frac 12}\|\phi\|_{L^{2}}
\end{equation*}
so that
\begin{equation*}
 \|S\|_{L^{2}\to L^{2}}\le
  |E|^{\frac 12-\frac 1p}|F|^{\frac 1q-\frac 12}
  g(T).
\end{equation*}
Since $g(T)\to0$, we see that for sufficiently large $T$ we have
$\|S\|_{L^{2}\to L^{2}}<1$, and the Amrein--Berthier estimate 
follows
as in the conclusion of the proof of Theorem \ref{the:DynAB1}.


\subsection{Proof of Theorem \ref{the:singPotNonMag}} 
\label{sub:singPotprf}

In order to adapt our method to the proof of Theorem 
\ref{the:singPotNonMag}, it is sufficient to note the following
facts.

\begin{itemize}
\item 
$S= \one{E} e^{ i T H} \one{F} e^{- i T H}$ is compact. 
Indeed, the operator $H$ is self-adjoint on $W^{2,2}$, 
$C_0^{\infty}$ is a core,
and the propagator $e^{i T H}$
can be represented as an integral operator with a bounded kernel 
$E(T,x,y)$, using e.g.~Theorem 1.1 in \cite{Yajima98-c}.
Thus $\one{E} e^{ i T H} \one{F}$ is an integral operator 
with an $L^2$ kernel, hence Hilbert-Schmidt, hence compact,
hence $S$ is a compact operator.

\item 
The operator $\cos(\sqrt{H+m^2})$ enjoys the finite speed of 
propagation property, thanks to \cite{Remling07-a},
since $H+m^2$ is essentially self-adjoint on 
$C_0^{\infty}(\mathbb{R}^n)$ and positive.

\item 
Unique continuation holds for $L^2$ solutions of
\begin{equation}
\label{eq:nonMagEig}
(H + m^2 )\phi = c \phi.
\end{equation}
To prove this, it is sufficient to invoke Theorem 6.3 of 
\cite{JerisonKenig85-a}, 
choosing as domain $\Omega$ a ball $B_R(0)$ with $R$ large enough.
We need only verify that $V \in L_{loc}^{\frac{n}{2}}$ and
$\phi \in W_{loc}^{2,\frac{2n}{n+2}}$. 
The second fact follows easily since 
$\phi \in L^2$ solves \eqref{eq:nonMagEig}, 
hence it is in the domain of $- \Delta + V$, that is
$ W^{2,2} \subset 
  W_{loc}^{2,\frac{2n}{n+2}}$.
It remains to check that $V\in L^{\frac{n}{2}}_{loc}$.
Since $1 < n^{\star} \leq 2$ and its H\"older conjugate exponent 
is $n-1$, condition \eqref{eq:V1Four} implies
\begin{equation*}
  \bra{x}^{1 + \sigma} \mathcal{F}^{-1} (\bra{\xi}^{\gamma} 
  \mathcal{F} V_{1,\epsilon}) \in L^{n-1} 
  \implies
  \mathcal{F}^{-1} (\bra{\xi}^{\gamma} 
  \mathcal{F} V_{1,\epsilon}) \in L^{n-1}
\end{equation*}
which in turn implies
\begin{equation*}
  V_{1,\epsilon} \in H^{\gamma,{n-1}} 
  \subset L^{n-1}
  \subset L_{loc}^{\frac{n}{2}}.
\end{equation*}
On the other hand,
$\mathcal{F} V_{2,\epsilon} $ is the difference of 
two positive measures $ \mu^{+}$ and $ \mu^{-}$. 
The total variation of $\mathcal{F} V_{2,\epsilon}$ is 
$(\mu^{+}+\mu^{-})(\mathbb{R}^n)$, and it is finite by assumption,
thus $\mu^{+}$, $\mu^{-}$ are finite measures. 
Since the (inverse) Fourier transform of a finite measure 
is continuous, 
$\mathcal{F}^{-1} \mu^{+} - \mathcal{F}^{-1} \mu^{-} = 
  V_{2,\epsilon} \in L_{loc}^{\frac{n}{2}}$. 
We conclude that 
$V = V_{1,\epsilon} + V_{2,\epsilon} \in 
  L_{loc}^{\frac{n}{2}}$.
\end{itemize}

\subsection{Proof of Theorem \ref{the:magtras}}\label{sub:magn}

Consider case (i).
Thanks to \eqref{eq:compat}, 
the vector field $\partial_v A$ is irrotational in $\mathbb{R}^n$,
hence conservative, and we can find
$\tilde{a} \in C^{\infty}(\mathbb{R}^n;\mathbb{R})$ such that
\begin{equation} \label{eq:magtraspot}
  \partial\tilde{a} = \partial_v A.
\end{equation}
Pick $\gamma\in C^{\infty}(\mathbb{R}^{n};\mathbb{R})$ such that
$\partial_{v}\gamma+\widetilde{a}=0$, for instance
\begin{equation*}
  \textstyle
  \gamma(x)=
  -\int^{(x \cdot v)}_{0} \tilde{a}(x + (s - (x \cdot v))v ) ds.
\end{equation*}
Then we have
\begin{equation*}
  0=\partial(\partial_{v}\gamma+\widetilde{a})=
  \partial_{v}(\partial \gamma+A)
\end{equation*}
so that $A+\partial \gamma$ is a constant vector over 
lines of direction $v$.
The gauge transform $u= e^{-i \gamma(x)}\tilde{u}$ gives
\begin{equation*}
  Hu=e^{-i \gamma(x)}\widetilde{H}\tilde{u}
  \quad\text{where}\quad 
  \widetilde{H}=(i \partial +A+\partial \gamma)^{2}+V.
\end{equation*}
We see that the coefficients of the modified operator are constant over lines of direction $v$. Since it is equivalent
to estimate $u$ or $\tilde{u}=e^{i \gamma(x)}u$, without loss of
generality we can assume that the
coefficients $A,V$ are constant over lines of direction $v$.

As in the proof of Theorem \ref{the:DynAB1}, it is sufficient to prove that the operator $S=\one{E} U(T,0) \one{F} U(0,T)$ on $L^{2}$ has norm strictly less than $1$, and again we can proceed by contradiction, assuming that $\lVert S \rVert = 1$. Here, $U(0,T)$ is the forward propagator of the magnetic Schr\"odinger equation, and $U(T,0)$ is the backward propagator.

Under our assumptions, there is $0 < \bar{T} \leq + \infty$ such that for any $0 < T < \bar{T}$ the backward propagator $U(T,0)$ can be represented as an integral operator with a bounded kernel $K(T,x,y)$ (see Theorem \ref{the:yaji}).

Now, for any finite measure sets $E$ and $F$, $\one{E} U(T,0) \one{F}$ is an integral operator with kernel $ \one{E}(x) K(T,x,y) \one{F}(y) \in L^2(\mathbb{R}^n \times \mathbb{R}^n)$. Thus it is an Hilbert-Schmidt operator, in particular it is compact, and as a consequence $S$ is compact as well. Repeating the steps in the proof of Theorem \ref{the:DynAB1}, we may construct an eigenfunction $g\in L^{2}$ for the operator $S$, corresponding to the eigenvalue 1.

Then we consider the translation operator
\begin{equation*}
  e^{i \tau ( i \partial_v)} \phi =  
  \phi(x - \tau v).
\end{equation*}
Since the coefficients are constant over lines of direction $v$, the translation
operator commutes with $H$ and hence for all $\tau,T$ we have
\begin{equation}
\label{eq:commPropTrans1}
[e^{i \tau ( i \partial_v)},U(0,T)]=0.
\end{equation}
\begin{equation}
\label{eq:commPropTrans2}
[e^{i \tau ( i \partial_v)},U(T,0)]=0.
\end{equation}
Indeed, if we take $f$ in the Schwartz space, $u = e^{i \tau ( i \partial_v)} U(0,t) f$ is a solution of  $i \partial_t u = H u$ with initial datum $e^{i \tau ( i \partial_v)} f$, hence (by uniqueness) it is equal to $U(0,t) e^{i \tau ( i \partial_v)} f$. By a density argument, \eqref{eq:commPropTrans1} holds. The argument for \eqref{eq:commPropTrans2} is analogous.

Define two sequences of sets
by translations in the direction $v$, as follows:
$E_{0}=E$, $F_{0}=F$, and for $k\ge1$
\begin{equation}\label{eq:defEk}
  E_{k+1}=E_{k}\cup(E_{k}+\tau_{k}v),
  \qquad
  F_{k+1}=F_{k}\cup(F_{k}+\tau_{k}v).
\end{equation}
Finally set $E_{\infty} = \bigcup_{k=0}^{\infty} E_{k}$, $F_{\infty} = \bigcup_{k=0}^{\infty} F_{k}$.
We may take a sequence $\{ \tau_{k} \}$ such that
\begin{equation*}
|E_k| < |E_{k+1}| < |E_k| + \frac{1}{2^{k}},
\qquad
|F_k| < |F_{k+1}| < |F_k| + \frac{1}{2^{k}},
\end{equation*}
so that $|E_{\infty}|$,  $|F_{\infty}| < + \infty$.

With this choice and thanks to properties 
\eqref{eq:commPropTrans1} and \eqref{eq:commPropTrans2}, 
$\{ e^{i \tau_l ( i \partial_v)} g \}_{l \in \mathbb{N}}$ 
is an infinite linearly independent set of eigenfunctions for 
$\one{E_{\infty}} U(T,0) \one{F_{\infty}} U(0,T) $ 
with the same eigenvalue 1, and this is in contradiction 
with the compactness of the operator.
This concludes the proof in case (i).

Consider now case (ii) in which $H$ is assumed to commute
with rotations. We proceed as in case (i), arguing by 
contradiction and constructing an eigenfunction $\phi\in L^{2}$
of the compact operator $S$ defined as above; in particular,
$\phi$ is supported on a subset $E'$ of $E$ and $e^{-iTH}$
in a subset $F'$ of $F$. Note that either $E'$ or $F'$
is \emph{not} rotationally invariant.

Given a matrix $M \in O(n)$, the  \emph{rotation} operator 
$\rho_{M}$ is the map $f \mapsto \rho_{M}(f)(x)=f(M x)$. 
If we consider a sequence of matrices $M_k \in O(n)$ 
such that $| M_k x - x| \rightarrow 0$ as 
$k \rightarrow + \infty$, it is easy to see that 
\begin{equation} \label{eq:L1convRot}
  \lVert \rho_{M_k} f - f\rVert_{L^1(\mathbb{R}^n)} \rightarrow 0 \text{ when } k \rightarrow + \infty, 
  \qquad \forall f \in L^1(\mathbb{R}^n).
\end{equation}
Given a finite measure set $A$, applying \eqref{eq:L1convRot} to $f=\one{A}$, up to the extraction of a subsequence from $M_k$, one has
\begin{align*}
  |M^{-1}_{k} A|  - | A| \leq 
    |M^{-1}_{k} A|   - |A \cap M^{-1}_{k} A| 	
  & = \int_{\mathbb{R}^n} \one{M^{-1}_{k} A \setminus A}(x) dx 
  \\ & \leq \int_{\mathbb{R}^n} \one{M^{-1}_{k} A \setminus A }(x) 
    + \one{A \setminus M^{-1}_{k} A }(x) dx 
  \\ & = \int_{\mathbb{R}^n} | (\rho_{M_k} \one{A})(x) 
    - \one{A }(x) | dx  \leq \frac{1}{2^{k}}.
\end{align*}
We now define two sequences of sets as follows:
$E_{0}=E$, $F_{0}=F$, and for $k\ge1$
\begin{equation}\label{eq:defEkb}
  E_{k+1}=E_{k}\cup(M_{k}^{-1} E_{k}),
  \qquad
  F_{k+1}=F_{k}\cup(M_{k}^{-1} F_{k}).
\end{equation}
Finally we set $E_{\infty} = \bigcup_{k=0}^{\infty} E_{k}$, $F_{\infty} = \bigcup_{k=0}^{\infty} F_{k}$.
We may take the sequence $M_k$ such that
\begin{equation*}
|E_k| < |E_{k+1}| < |E_k| + \frac{1}{2^{k}},
\qquad
|F_k| < |F_{k+1}| < |F_k| + \frac{1}{2^{k}},
\end{equation*}
so that $|E_{\infty}|$,  $|F_{\infty}| < + \infty$.

With this choice, consider 
$\{ \rho_{M_k} \phi \}_{k \in \mathbb{N}}$. 
This is a linearly independent set, since the support of 
$\phi$ is not rotationally invariant.
Moreover, thanks to the fact that $U(0,T)$ and $\rho_{M_k}$ 
commute, each element of the set is an eigenfunction for 
$\one{E_{\infty}} U(T,0) \one{F_{\infty}} U(0,T)$ 
with the same eigenvalue 1, 
and this is in contradiction with the compactness of the operator,
proving the Theorem in case (ii).

\subsection{Proof of Theorem \ref{the:Obs1}}\label{sub:Obs1}

In order to prove \eqref{eq:Obs2}
it is sufficient to prove the observability inequality
\begin{equation} \label{eq:Obs}
  \lVert u_{0} \rVert_{L^2} \leq C \int_0^T 
  \lVert e^{- i t H} u_{0} \rVert_{L^2(E^c)} dt
\end{equation}
and then apply the H\"{o}lder inequality in time.

Write $f(t) = \lVert e^{- i t H} u_{0} \rVert_{L^2(E^c)}$. 
By Theorem \ref{the:DynAB1} 
(resp. Theorem \ref{the:singPotNonMag}), 
we know that there can be at most one time 
$t^{\star}(T) \in [0,T]$ such that $f(t^{\star}(T))=0$.
Thus we can certainly find four times 
\begin{equation*}
  0 \leq t_{1}(T)<t_{2}(T)<t_{3}(T)<t_{4}(T) \leq T 
\end{equation*}
such that 
$\min_{t \in [t_{1}(T),t_{2}(T)]} f(t) = f(\tau_{1}(T)) >0$, 
$\min_{t \in [t_{3}(T),t_{4}(T)]} f(t) = f(\tau_2(T)) >0$, 
where $\tau_{1}(T) \in [t_{1}(T),t_{2}(T)]$, 
$\tau_{2}(T) \in [t_{3}(T),t_{4}(T)]$ are the points of minimum. 
Hence
\begin{align*} 
  \textstyle
  \int_{0}^T f(t) dt \geq  &  
  \textstyle
  \int_{t_{1}(T)}^{t_{2}(T)} f(t)dt 
  + \int_{t_{3}(T)}^{t_{4}(T)} f(t) dt 
\\ 
  & \geq [t_{2}(T) - t_{1}(T)] f(\tau_1(T)) 
  + [t_{4}(T) - t_{3}(T)]f(\tau_2(T)) 
\\
  &
  \geq \min \{t_{2}(T) - t_{1}(T), t_{4}(T) 
  - t_{3}(T) \} (f(\tau_1(T)) + f(\tau_2(T))).
\end{align*}
Recalling the definition of $f$ and applying \eqref{eq:AmBerIneq}
we obtain
\begin{equation*}
\int_0^T \lVert e^{-i t H} u_{0} \rVert_{L^2(E^c)} dt \geq 
\frac{ \min \{t_{2}(T) - t_{1}(T), t_{4}(T) - t_{3}(T)) \}}{C(E,A,V,\tau_2(T)-\tau_1(T))} \lVert e^{- i \tau_{1}(T) H} u_{0} \rVert_{L^2},
\end{equation*}
and \eqref{eq:Obs} follows easily.

\section{Metaplectic operators}\label{sec:metOp}

We recall some basic definitions.
A \emph{symplectic matrix} is a matrix
$S \in \mathbb{R}^{2n \times 2n}$ 
such that $S^{t} J S = J$, where
\begin{equation} \label{eq:JskewSym}
  J=
  \begin{pmatrix}
    0 & I \\ 
    -I & 0
  \end{pmatrix}.
\end{equation}
If we decompose $S$ in $n \times n$ blocks
\begin{equation} \label{eq:sympMat}
  \begin{pmatrix}
    A & B \\ 
    C & D
  \end{pmatrix},
\end{equation}
an equivalent condition is
\begin{equation}
\label{eq:sympCond}
CD^{t} = DC^t, \qquad AB^t = BA^t, \qquad AD^t - BC^t = I.
\end{equation}
When the block $B$ is invertible, 
the matrix is said to be \textit{free}.

Now, consider the phase space coordinates 
$(x,\xi) \in \mathbb{R}^{2n}$ and $\tau \in \mathbb{R}$. 
The \emph{Heisenberg group} can be defined as $\mathbb{R}^{2n+1}$
endowed with the group law
\begin{equation*}
  \textstyle
  (x,\xi;\tau)\cdot(x',\xi';\tau')=
  (x+x',\xi+\xi';\tau+\tau'+\frac 12(x \cdot \xi'-\xi \cdot x')).
\end{equation*}
The Heisenberg group can be represented as a group 
of bounded operators on $L^{2}(\mathbb{R}^n)$ via the
\emph{Schr\"{o}dinger representation}
$(x,\xi;\tau)\mapsto \rho(x,\xi;\tau)$ defined as
\begin{equation} \label{eq:heisRep}
  \rho(x,\xi;\tau) \phi(y) = 
  e^{2 \pi i \tau} e^{- i \pi x \cdot \xi} 
  e^{2 \pi i\xi \cdot y} \phi(y - x).
\end{equation}
Now, given a symplectic matrix $S\in \mathbb{R}^{2n \times 2n}$,
the metaplectic operator $\hat{S}$ with projection 
$\pi(\hat{S}) = S$ is the unitary operator, 
unique up to a constant phase factor, such that
\begin{equation} \label{eq:metDef}
  \hat{S} \rho(x,\xi;\tau) \hat{S}^{-1} = \rho(S(x,\xi),\tau).
\end{equation}
Given a symplectic matrix $S$, we will denote by $\mu(S)$ 
the corresponding metaplectic operator, 
defined through \eqref{eq:metDef}.

We shall prove the following:

\begin{theorem}[Metaplectic Amrein-Berthier inequality]
  \label{the:metAB}
  Let $S$ be a free symplectic matrix and $\hat{S}$ a metaplectic operator with $\pi(\hat{S})=S$. Then for each couple of sets of finite measure $E$, $F$ the following inequality holds:
  \begin{equation} \label{eq:metOpABIneq}
    \|\phi \|_{L^{2}}\le
          C (  \|\one{E^{c}} \phi \|_{L^{2}} 
          +  \|\one{F^{c}} \hat{S} \phi \|_{L^{2}} )
        \qquad
        \forall t\in \mathbb{R}.
  \end{equation}
\end{theorem}

In particular, if both $f$, $\hat{S} f$ are supported 
on sets of finite measure, then $f=\hat{S}f=0$ 
(see Theorem 6.2 in \cite{GrochenigShafkulovska24-a}).

The proof of Theorem \ref{the:metAB} follows the usual scheme, 
using a few supporting lemmas.

\begin{lemma} \label{lem:metOpComp}
  Under the assumptions of Theorem \ref{the:metAB} 
  the operator $\one{E} \hat{S}^{-1} \one{F} \hat{S}$ is compact.
\end{lemma}

If the projection of the metaplectic operator is free, 
the operator  admits the following representation formula 
(see Remark 3.3 in \cite{CorderoGiacchiMalinnikova24-a}):

\begin{proposition}[Representation Formula for 
  free symplectic matrices] \label{prop:freeRep}
  Given a metaplectic operator $\hat{S}$ with free projection $S$, 
  for every $\phi \in L^2(\mathbb{R}^n)$ we have:
  \begin{equation} \label{eq:freeRep}
    \hat{S} \phi (y) = 
    | \operatorname{det}(B) |^{-\frac{1}{2}} 
    e^{i \pi D B^{-1} y \cdot y} \int_{\mathbb{R}^n} \phi(t) 
    e^{i \pi B^{-1} A t \cdot t} 
    e^{- 2 \pi i (B^{-1} y \cdot t )} dt ,
  \end{equation}
  with respect to the block decomposition \eqref{eq:sympMat}.
\end{proposition}

\begin{proof}[Proof of Lemma \ref{lem:metOpComp}]
  It is sufficient to prove the compactness of 
  $\one{E} \hat{S}^{-1} \one{F}$. 
  Now, if $S^{-1}$ is free, 
  one can just plug the integral representation formula 
  given in Proposition \ref{prop:freeRep} 
  to find out that 
  $\one{E} \hat{S}^{-1} \one{F}$
  is an integral operator 
  with a kernel in $L^2(\mathbb{R}^{n} \times \mathbb{R}^n)$, 
  and thus it is Hilbert-Schmidt, and in particular compact.
\end{proof}

\begin{remark} \label{rem:metCompFail}
  If the block $B \neq 0$ is not invertible, 
  compactness of the metaplectic operator can fail. 
  Consider the following easy counterexample (in dimension $2$, 
  but it can be easily generalized). 
  Choose $E=F=Q$ a square: 
  \begin{equation*}
  \one{Q}(x)= \one{[0,1]}(x_1) \one{[0,1]}(x_2)
  \end{equation*}
  and $\hat{S}$ with projection $S$ with blocks:
  \begin{equation*}
  A = D =
  \begin{pmatrix}
  1 & 0\\
  0 & 0
  \end{pmatrix},
  \qquad
  B = - C =
  \begin{pmatrix}
  0 & 0\\
  0 & 1
  \end{pmatrix}
  \end{equation*}
  i.e., $\hat{S} = \mathcal{F}_2$ is the (partial) 
  Fourier transform in the second variable.
  Note that $L^2(\mathbb{R}^2;\mathbb{C})$ 
  can be identified with the tensor product 
  $L^2(\mathbb{R};\mathbb{C}) \bigotimes L^2(\mathbb{R};\mathbb{C})$. 
  Elementary tensors are:
  \begin{equation*}
    f_1 \otimes f_2 = f_1(x_1)f_2(x_2)
  \end{equation*}
  The operator $\chi_Q \hat{S} \chi_Q$ itself tensorizes, i.e:
  \begin{equation*}
    \one{Q} \hat{S} \one{Q} (f_1 \otimes f_2)= 
    \one{[0,1]} f_1(x_1) \otimes 
      (\one{[0,1]} \mathcal{F} \one{[0,1]} f_2)(x_2)
  \end{equation*}
  We can thus write:
  \begin{equation*}
    \one{Q} \hat{S} \one{Q} = 
    \one{[0,1]} \otimes \one{[0,1]} \mathcal{F} \one{[0,1]}
  \end{equation*}
  By Theorem 2 in \cite{ZanniKubrusly15-a}, 
  if the operator were compact, then $\one{[0,1]}$ 
  would be a compact operator, but clearly it is not one.
\end{remark}

\begin{remark} \label{rem:ineqFail}
  Under the assumptions of Remark \ref{rem:metCompFail}, 
  one can say something more, and actually construct 
  counterexamples to inequality \eqref{eq:metOpABIneq}. 
  Consider the partial Fourier transform $
  \hat{S}= \mathcal{F}_2$ in the second variable in dimension $2$,
  introduced in Remark \ref{rem:metCompFail}. Define
  \begin{equation} \label{eq:countexSet}
    E^{\alpha}_{a} = \{ (x_1,x_2) \in \mathbb{R}^2 | 0 
    \leq x_1 \leq a , -x_1^{\alpha} \leq x_2 \leq x_1^{\alpha} \}
  \end{equation}
  where $-1 < \alpha < 0$ is fixed and $0 < a \leq 1$ 
  is a parameter, and
  \begin{equation} \label{eq:countexFunc}
    f_{a,\alpha}(x_1,x_2)=\one{E^{\alpha}_{a}}(x_1,x_2) 
    = \one{[0,a]}(x_1) \one{[-x_1^{\alpha},x_1^{\alpha}]}(x_2)
  \end{equation}
  Since $x_1^{\alpha}$ is integrable, 
  the set $E^{\alpha}_{a}$ has finite measure. Then we have:
  \begin{equation} \label{eq:ineqFailnorm1}
    \textstyle
    \lVert f_{a,\alpha} \lVert_{L^2(\mathbb{R}^n)} = 
    | E^{\alpha}_{a}|^{\frac{1}{2}} = 
    ( 2 \int_{0}^{a} x^{\alpha} )^{\frac{1}{2}} 
    = \sqrt{\frac{2}{\alpha+1}}\cdot a^{\frac{\alpha+1}{2}}.
\end{equation} 

Consider the finite measure set $E=E^{\alpha}_{1}$. 
Since $E_{a}^{\alpha} \subset E$, $f_{a,\alpha}$ 
vanishes on $E^{c}$, so that
\begin{equation} \label{eq:ineqFailnorm2}
  \lVert f_{a,\alpha} \rVert_{L^2(E^{c})} = 0
\end{equation}
On the other hand, we have easily
\begin{equation*}
  \textstyle
\mathcal{F}_2 f_{a,\alpha}(x_1,x_2) = \one{[0,a]}(x_1) \frac{\sin(2 \pi x_1^{\alpha} x_2)}{\pi x_2}.
\end{equation*}
Since $\mathcal{F}_2 f_{a,\alpha}(x_1,x_2)$ 
vanishes outside the strip $0 \leq x_1 \leq a$, we can write
\begin{equation}\label{eq:ineqFailnorm3}
  \lVert \mathcal{F}_2 f_{a,\alpha} (x_1,x_2) \rVert_{L^2(E^c)} = 
  \textstyle
  \lVert \one{[0,a]} (x_1) \frac{\sin(2 \pi x_1^{\alpha} x_2)}{\pi x_2} \rVert _{L^2(A)}
  =\sqrt{2} (\int_{0}^{a} \int_{x_1^{\alpha}}^{+ \infty} \frac{\sin^2(2 \pi x_1^{\alpha} x_2)}{\pi^2 x_2^2} dx_2 dx_1 )^{\frac{1}{2}} 
\end{equation}
where $A=\{x\in \mathbb{R}^{2}:0\le x_{1}\le a,\ 
  x_2 \not\in [-x_1^{\alpha} , x_1^{\alpha}]\}$.
Setting $y = x_1^{\alpha} x_2$ we get
\begin{equation*}
  \textstyle
=\frac{\sqrt{2}}{\pi} ( \int_{0}^{a} x_1^{\alpha} \int_{x_1^{2 \alpha}}^{+ \infty} \frac{\sin^2(2 \pi y)}{y^2} dy dx_1 )^{\frac{1}{2}} \leq \frac{\sqrt{2}}{\pi}(\int_{0}^{a} x_1^{\alpha} \int_{x_1^{2 \alpha}}^{+ \infty} y^{-2} dy dx_1 )^{\frac{1}{2}} =\frac{\sqrt{2}}{\pi} (\int_{0}^{a} x_1^{- \alpha} dx_1 )^{\frac{1}{2}}
\end{equation*}
and in conclusion
\begin{equation*}
\lVert \mathcal{F}_2 f_{a,\alpha} (x_1,x_2) \rVert_{L^2(E^c)} \le \frac{\sqrt{2}}{\pi \sqrt{(- \alpha + 1)}} a^{\frac{-\alpha + 1}{2}}.
\end{equation*}
If the inequality \eqref{eq:metOpABIneq} were true 
with a constant $C$ independent on $a$, 
plugging $f_{a,\alpha}$ in the inequality would give
\begin{equation*}
  a^{\frac{\alpha + 1}{2}} \leq 
  \frac{C}{\pi} \sqrt{\frac{\alpha + 1}{-\alpha+1}} 
  a^{\frac{- \alpha + 1}{2}} \implies a^{- |\alpha|} 
  \leq \frac{C}{\pi} \sqrt{\frac{-|\alpha| + 1}{|\alpha| + 1}}
\end{equation*}
which is clearly false for $a$ small.
\end{remark}

\begin{lemma} \label{lem:eigfnsMet}
  Assume $\hat{S}$, $E$, $F$ are as in Theorem \ref{the:metAB},
  and let $g$ be a $1$-eigenfunction of the operator 
  $\one{E} \hat{S}^{-1} \one{F} \hat{S}$.
  Then there are a sequence 
  $\{ y_k \}_{k \in \mathbb{N}}\subset \mathbb{R}^n$, 
  and two finite measure sets $E_{\infty} \supset E$ 
  and $F_{\infty} \supset F$ such that 
  $\{ g(\boldsymbol{\cdot} - y_k) \}$ is a family of 
  linearly independent $1$-eigenfunctions of the operator 
  $\one{E_{\infty}} \hat{S}^{-1} \one{F_{\infty}} \hat{S}$.
\end{lemma}

\begin{proof}[Proof of Lemma \ref{lem:eigfnsMet}]
  Similarly to the proof of Theorem \ref{the:magtras},
  define the sets 
  $E_{0}=E$, $E_{k+1}=E_{k}\cup(E_{k}+y_{k})$
  and
  $F_{0}=F$, $F_{k+1}=F_{k}\cup(F_{k}+Ay_{k})$,
  where $A$ is the top-left block of $S$ 
  with respect to the block decomposition \eqref{eq:sympMat}. 
  We can choose the sequence $\{ y_{k} \}_{k \in \mathbb{N}}$ 
  with the property
  \begin{equation*}
  |E_k| < |E_{k+1}| < |E_k| + \frac{1}{2^{k}},
  \qquad
  |F_k| \leq |F_{k+1}| \leq |F_k| + \frac{1}{2^{k}},
  \end{equation*}
  so that $|E_{\infty}|$,  $|F_{\infty}| < + \infty$
  Invoking definition \eqref{eq:metDef}, we see that
  \begin{align*} 
    ( \one{E_{\infty}} \hat{S}^{-1} \one{F_{\infty}} \hat{S}) (g(\boldsymbol{\cdot} - y_{k})) &= \one{E_{\infty}} \hat{S}^{-1} \one{F_{\infty}} \hat{S} \rho(y_{k},0;0) g
  \\ 
    & = \one{E_{\infty}} \hat{S}^{-1} 
    \one{F_{\infty}} \rho(S(y_{k},0);0) \hat{S} g 
  \\
    &= \one{E_{\infty}} \hat{S}^{-1} \rho(S(y_{k},0);0) \one{F_{\infty}} \hat{S} g
  \\
    & = \one{E_{\infty}} \rho(y_{k},0;0) \hat{S}^{-1} 
    \one{F_{\infty}} \hat{S} g
  \\ 
  & = \rho(y_{k},0;0) \one{E_{\infty}} \hat{S}^{-1} \one{F_{\infty}} \hat{S} g
  \\ 
  &= \rho(y_{k},0;0) g = g (\boldsymbol{\cdot} - y_{k}).
  \qedhere
  \end{align*}
\end{proof}

\begin{proof}[Conclusion of the proof]
  Combining Lemma \ref{lem:metOpComp} and Lemma 
  \ref{lem:eigfnsMet}, and proceeding as in the proof of 
  Theorem \ref{the:DynAB1}, we conclude the proof of
  Theorem \ref{the:metAB}.
\end{proof}

We can apply the last results in the dynamical setting. 
Indeed, given a symmetric real $2n \times 2n$ matrix $M$, 
we associate with $M$
the classical hamiltonian
\begin{equation} \label{eq:clasHam}
  H_{C} = \frac{1}{2} 
  \langle M (x,\xi) , (x,\xi)\rangle_{\mathbb{R}^{2n}}.
\end{equation} 
and we can naturally define a quantum hamiltonian 
by means of the Weyl quantization:
\begin{equation} \label{eq:WeylHam}
  H f (x) = \int_{\mathbb{R}^{2n}} e^{2 \pi i (x-y) \cdot \xi} 
  H_{C}(\textstyle\frac{x+y}{2} , \xi) f(y) dy d \xi
\end{equation}
One can show (see Section 5 in 
\cite{CorderoGiacchiMalinnikova24-a}) 
that the propagator $e^{- i T H}$ at any time $T$ is metaplectic, 
and its symplectic projection is the matrix exponential of $T J M$:
\begin{equation*}
\pi(e^{- i T H}) = e^{T J M}.
\end{equation*}

\begin{corollary} \label{cor:dynQuad}
  Let $H$ be as in \eqref{eq:WeylHam}, and let
  $u(t,x) = e^{- i t H} f$ be the unique solution to 
  $i \partial_t u = H u$ with
  $u(0) = f\in L^{2}(\mathbb{R}^n)$.
  Then for all $T>0$ the following holds:
  \begin{itemize}
  \item 
  If the top right block $B$ (with respect 
  to the block decomposition \eqref{eq:sympMat}) of $e^{TJM}$
   is $\neq 0$, then $u(0)$ and $u(T)$ cannot 
   be both supported over a finite measure set.
  \item 
  If in addition the block $B$ is invertible, 
  the inequality \eqref{eq:AmBerIneq} holds 
  for any sets $E$ and $F$ of finite measure.
  \end{itemize}
\end{corollary}

\begin{proof}
  The claim follows from Theorem 6.2 of
  \cite{GrochenigShafkulovska24-a} when $B \neq 0$ 
  and from Theorem \ref{the:metAB} in this paper 
  when $B$ is invertible. 
\end{proof}

We can give the following applications to 
generic quantum harmonic oscillators 
(see again Section 5 of \cite{CorderoGiacchiMalinnikova24-a}).

\begin{example}[Harmonic Oscillator]
\label{ex:genHarmOsc}
Consider the Schr\"odinger equation 
for the quantum harmonic oscillator:
\begin{equation}
  \textstyle
  i \partial_t u(t,x) = 
  ( - \Delta + \sum_{j=1}^n (\frac{\omega_j}{2})^{2} x^2_j ) u(t,x)
\end{equation}
for $\omega_{j} \neq 0$, $j = 1 , \ldots , n$.
Then one can verify that
\begin{equation*}
  \pi(e^{- i T H}) =
  \begin{pmatrix}
    A_T & B_T \\
    C_T & D_T
  \end{pmatrix}
\end{equation*}
where
\begin{equation*}
  A_T = D_T = 
  \operatorname{diag} (\{ \cos(\omega_j T) \}_{j=1,\ldots,n} )
\end{equation*}
\begin{equation*}
  B_T =  \operatorname{diag} 
  (\{ \frac{ \sin(\omega_j T)}{\omega_j} \}_{j=1,\ldots,n}).
\end{equation*}
\begin{equation*}
  C_T = 
 -  \operatorname{diag} 
  (\{ \omega_j \sin(\omega_j T)\}_{j=1,\ldots,n}).
\end{equation*}
In particular, $B_T \neq 0$ if 
$\exists j \in \{ 1, \ldots , n \}$ such that 
$\frac{\omega_j T}{\pi}$ is not an integer.
In this case, 
the first conclusion of Corollary \ref{cor:dynQuad} holds. 
Note that if the oscillator is \textit{anisotropic}, 
i.e. not all the frequencies are integer multiples 
of the same number, the uncertainty principle 
(without the inequality) holds for any $T > 0$.

Moreover, $B$ is invertible if $\frac{\omega_j T}{\pi}$ 
is not an integer for any $j = 1 , \ldots , n$. 
In this case, the second conclusion of Corollary \ref{cor:dynQuad} 
holds and we recover inequality \eqref{eq:AmBerIneq}. 
In particular, the set of times $T$ 
for which the inequality can fail is discrete.
\end{example}

\begin{example}[Partial Harmonic Oscillator] \label{ex:partHarmOsc}
For simplicity, we consider the case $n=2$ and the equation
\begin{equation}
  \textstyle
  i \partial_t u(t,x) = ( - \frac{1}{4 \pi}\partial^2_{x_2}
   + \pi x_2^2) u(t,x).
\end{equation}
Defining
\begin{equation*}
A_T = D_T =
\begin{pmatrix}
1 & 0 \\
0 & \cos(T)
\end{pmatrix},
\end{equation*}
\begin{equation*}
B_T = - C_T =
\begin{pmatrix}
0 & 0 \\
0 & \sin(T)
\end{pmatrix},
\end{equation*}
one can verify that
\begin{equation*}
\pi(e^{- i T H}) =
\begin{pmatrix}
A_T & B_T \\
C_T & D_T
\end{pmatrix}.
\end{equation*}
We see that for times $T \neq k \pi$, $k$ integer, 
the block $B_T$ is $\neq 0$, 
and thus the first conclusion of Corollary \ref{cor:dynQuad} holds.
\end{example}


\end{document}